\documentclass[english]{article} \usepackage{a4wide}
\usepackage{lmodern} \usepackage[T1]{fontenc}
\usepackage[latin9]{inputenc} 
\usepackage{color} \usepackage{graphicx} \usepackage{amssymb}
\usepackage{xspace} \usepackage{hyperref} \usepackage{bbm}
\usepackage{babel} \usepackage{oupbib}

\newcommand{\ie}{{\em i.e.\/}\xspace}
\newcommand{\eg}{{\em e.g.\/}\xspace}

\newcommand{\etc}{{\em etc.\/}\xspace}

\newcommand{\todo}[1]{{\bf  *** #1 ***}}

\newcommand{\HIDE}[1]{ }
\newcommand{\COMMENT}[1]{ }

\newcommand{\half}{\frac{1}{2}}

\newcommand{\eqref}[1]{\mbox{(\ref{eq:#1})}}

\newcommand{\secref}[1]{\mbox{\S$\,$\ref{sec:#1}}}

\newcommand{\itref}[1]{\mbox{\ref{it:#1}}}

\newcommand{\condref}[1]{\mbox{Condition~\ref{cond:#1}}}

\newcommand{\nei}[1]{{\rm ne}(#1)}

\newcommand{\bd}{{\rm bd}}

\newcommand{\transp}{^{\prime}}

\newcommand{\cd}{\,|\,}

\newcommand{\E}{{\mbox{E}}}

\newcommand{\var}{{\mbox{var}}}

\newcommand{\cip}{\mbox{$\perp\!\!\!\perp$}}

\newcommand{\R}{{\cal R}}

\newtheorem{expl}{Example}
\newtheorem{definer}{Definition}
\newtheorem{theorem}{Theorem}
\newtheorem{thm}[theorem]{Theorem}
\newtheorem{algor}{Algorithm}

\newtheorem{cond}{Condition}
\newtheorem{rem*}{Remark}

\newcommand{\halm}{\hspace*{\fill} $\Box$\par}

\newenvironment{ex}{\begin{expl}\rm}{\halm\end{expl}}

\newenvironment{defn}{\begin{definer}\rm}{\end{definer}}

\newcommand{\iid}{{independent and identically distributed}\xspace}

\renewcommand{\todo}[1]{}
\newcommand{\reals}{{\mathbbm R}}

 \newcommand{\bmu}{{\bm
    \mu}} \newcommand{\cX}{{\cal X}}

\newcommand{\bY}{{\bm Y}} \newcommand{\by}{{\bm y}}
\newcommand{\bx}{{\bm x}} \renewcommand{\bd}{{\bm d}}
\newcommand{\bo}{{\bf 0}} \newcommand{\norm}{{\cal N}}
\newcommand{\BF}{\mbox{\rm BF}} \newcommand{\rss}{\mbox{\rm
    RSS}\xspace} \newcommand{\aic}{\mbox{\rm AIC}\xspace}
\newcommand{\bm}[1]{\mbox{\boldmath $#1$}} \newcommand{\btheta}{{\bm
    \theta}} \newcommand{\bbl}{\mbox{$\mathbbm 1$}}
\newcommand{\bX}{{\bm X}} \newcommand{\bt}{{\bm t}}
\newcommand{\bT}{{\bm T}} \newcommand{\bnabla}{{\bm \nabla}}
\newcommand{\IF}{\mbox{\rm IF}\xspace} \newcommand{\SF}{\mbox{\rm
    SF}\xspace} \newcommand{\SD}{\mbox{\rm SD}\xspace}
\newcommand{\tr}{\mathop{\mathrm{tr}}} \renewcommand{\R}{\reals}
\title{Theory and Applications of Proper Scoring Rules}

\author{{A. Philip Dawid, University of Cambridge}\\and\\{Monica
    Musio, Universit\`a di Cagliari}}

\date{}

\begin{document}
\maketitle

\begin{abstract}
  We give an overview of some uses of proper scoring rules in
  statistical inference, including frequentist estimation theory and
  Bayesian model selection with improper priors.
\end{abstract}

\section{Introduction}
\label{sec:intro}

The theory of {\em proper scoring rules\/} (PSRs) originated as an
approach to probability forecasting \cite{apd:enc(probfore)}, the
general enterprise of forming and evaluating probabilistic
predictions.  The fundamental idea is to develop ways of motivating a
forecaster to be honest in the predictions he announces, and of
assessing the performance of announced probabilities in the light of
the outcomes that eventuate.  The early applications of this
methodology were to meteorology \cite{brier:50} and subjective
Bayesianism \cite{good:52,definetti:75}.  However, it is becoming
clear that the mathematical theory of PSRs has a wide range of other
applications in Statistics generally.  The aim of this paper is to
give a brief outline of some of these.  Some further details will
appear in \textcite{DMV:2014}.

After setting out the basic ideas and properties of PSRs in
\secref{psr}, in \secref{special} we describe a number of important
special cases.  Section~\ref{sec:est} describes how a PSR supplies an
alternative to the likelihood function of a statistical model, and how
this can be used to develop new estimators that may have useful
properties of robustness and/or computational feasability.  In
particular, in \secref{evading} we show that, by using an appropriate
PSR, it is possible to avoid difficulties associated with intractable
normalising constants.  In \secref{modsel} the same ability to ignore
a normalising constant is shown to supply an approach to Bayesian
model selection with improper priors.  Some concluding comments are
made in \secref{conc}.

\section{Proper scoring rules}
\label{sec:psr} The most natural context in which to introduce the
idea of a proper scoring rule is that of a game between a
decision-maker (henceforth ``You'') and Nature.  Let $X$ be a random
variable with values in $\mathcal{X}$, and let ${\cal P}$ be a family
of distributions over ${\cal X}$.  In due course, Nature will reveal
the value, $x$, of $X$.  Ahead of that, Your task is to quote a
distribution $Q\in{\cal P}$, intended to represent Your uncertainty
about how $X$ might turn out.  Later, after Nature has revealed $x$,
You will suffer a penalty $S(x,Q)$, depending on both Your quoted
distribution $Q$, and Nature's chosen value $x$, for $X$.  The
function $S$ is a {\em scoring rule\/}.

Suppose Your actual beliefs about $X$ are represented by the
distribution $P\in{\cal P}$.  If Your quoted distribution is $Q$, You
assess Your expected penalty as
\begin{equation}
  \label{eq:spq}
  S(P,Q) := \E_{X\sim P}S(X,Q).
\end{equation}
According to the principles of decision theory, You should choose Your
quote $Q$ to minimise Your expected score $S(P,Q)$.  However, this
might or might not coincide with Your true belief $P$.  A {\em proper
  scoring rule (PSR)\/} is one that encourages You to be honest:
\begin{defn}
  The scoring rule $S$ is {\em proper\/} with respect to ${\cal P}$
  if, for $P,Q\in\mathcal{P}$, the expected score $S(P,Q)$ is
  minimised in $Q$ at $Q=P$.  Further $S$ is {\em strictly proper\/}
  if this is the unique minimum: $S(P,Q)>S(P,P)$ for $Q\neq P$.
\end{defn}

\subsection{General construction}
\label{gen} We have introduced proper scoring rules in the context of
a special kind of decision problem, where Your decision has the form
of a quoted distribution $Q$.  But virtually any decision problem
gives rise to an associated PSR.

Thus consider a decision problem, with state space ${\cal X}$ and
arbitrary action space ${\cal A}$.  Again Your task is to choose a
decision, this time of the form $a\in{\cal A}$, after which Nature
will reveal the value $x$ of $X$, and You will be subject to a loss
$L(x,a)$.  The loss function $L$ is, again, essentially arbitrary.

Let ${\cal P}$ be a family of distributions over ${\cal X}$ such that
for each $P\in{\cal P}$ there exists a {\em Bayes act\/}:
$$a_P := \arg\min_{a\in{\cal A}}L(P,a)$$
where $L(P,a) := \E_{X\sim P} L(X,a)$.  If the Bayes act is not
unique, we arbitrarily nominate one such act as $a_P$.  Now define a
scoring rule $S$ by:
\begin{equation}
  \label{eq:sdef}
  S(x,Q) = L(x, a_Q)\quad\quad(x\in{\cal X}, Q\in{\cal P}).
\end{equation}
Then $S(P,Q) = L(P, a_Q) \geq L(P, a_P) = S(P,P)$.  Thus $S$ is a
PSR with respect to ${\cal P}$.

In this way the theory of proper scoring rules subsumes a large part
of statistical decision theory.

\subsection{Related concepts}
\label{sec:related}

Let $S$ be a PSR with respect to a large convex family
${\cal P}$ over ${\cal X}$.  Starting from $S$, we can define a
collection of useful statistical functions:

\begin{description}
\item[Entropy] The minimised value $H(P) := S(P,P)$ is the {\em
    (generalised) entropy\/} of $P\in{\cal P}$.
\item[Discrepancy] The excess score $D(P,Q) := S(P,Q) - H(P)$ is the
  {\em discrepancy\/} or {\em divergence} of $Q\in{\cal P}$ from
  $P\in{\cal P}$.
\item[Metric] Locally, $D(P, P+dP)$ defines a {\em Riemannian
    Metric\/} on $\mathcal{P}$.
\item[Dependence Function] The dependence of $X$ on a random variable
  $U$ jointly distributed with $X$ is
    $$C(X,U) := H(P_X) - \E_U\{H(P_{X\cd U})\}.$$
  \end{description}

  We note without proof the following properties of these associated
  functions \cite{rr139b}:
  \begin{thm}\quad\\[-1ex]
    \begin{enumerate}
    \item $H(P)$ is a concave functional of $P$ (strictly concave if
      $S$ is strictly proper).
    \item $D(P,Q) \geq 0$, and $D(P,Q)-D(P,Q_0)$ is an affine function
      of $P$.
    \item For a parametric model
      $\{P_\btheta:\btheta\in\Theta\subseteq \R^p\}$, the metric takes
      the form
$$D(P_\btheta, P_{\btheta+{\rm d}\btheta}) = {\rm d}\btheta\transp g(\btheta) {\rm d}\btheta,$$
where the $(p \times p)$ matrix $g(\btheta)$ satisfies
 $$\partial g_{ab}(\btheta)/\partial \theta_c = \partial g_{ac}(\btheta)/\partial \theta_b.$$
\item $C(X,U)\geq 0$ and vanishes if $X\, \cip \, U$ (and, when $S$ is
  strictly proper, only in this case.)
\end{enumerate}
\end{thm}

We could alternatively start with an entropy function $H$, a
discrepancy function $D$, or a metric $g$, having the additional
properties described above.  In each case we can (under some technical
conditions) construct a PSR $S$ from which it can be derived
\cite{rr139b}.

\section{Some special proper scoring rules}
\label{sec:special}

Since every decision problem induces a PSR, there is a very great
variety of these: the set of PSRs has essentially the same cardinality
as the set of concave functionals (serving as associated generalised
entropy functions) on ${\cal P}$.  Here we discuss some cases of
special interest.  For further special cases, see among others
\textcite{rr139b,apd:psr,apd/ps:design}.  Where appropriate, we equip
${\cal X}$ with an underlying measure $\mu$ dominating ${\cal P}$, and
write $p(\cdot)$ for the density (Radon-Nikodym derivative) $dP/d\mu$,
\etc

\subsection{Log score}
\label{sec:logscore}

The {\em log score\/} \cite{good:52} is just negative log likelihood:

\begin{equation}
  \label{eq:logscore}
  S(x,Q) = -\ln q(x).
\end{equation}

For this case we find:
\begin{itemize}
\item $H(P) = - \int\! d\mu(y)\cdot p(y) \ln p(y)$ is the {\em Shannon
    Entropy\/} of $P$.
\item $D(P,Q) = \int\! d\mu(y)\cdot p(y) \ln \{p(y)/q(y)\}$ is the {\em
    Kullback-Leibler Discrepancy\/} $K(P,Q)$.
\item $C(X,U) = \E\left[\ln \left\{p(X,U)/p(X)p(U)\right\}\right]$ is
  the {\em mutual Information\/} $I(X;U)$.
\item $g(\btheta)$ is the {\em Fisher Information\/} matrix.
\end{itemize}

The well-known property that $K(P,Q) \geq 0$, with equality if and
only if $Q=P$, shows that the log score is strictly proper.

It is interesting to see so many fundamental ingredients of
statistical theory and information theory flowing naturally from the
log score.  But it is equally true that many of the important
properties for which these are renowned remain valid for the more
general constructions of \secref{related}.  In particular, the whole
theory of {\em information geometry\/}, which subsumes but goes beyond
the information metric on ${\cal P}$, can be generalised to yield the
{\em decision geometry\/} associated with a given PSR
\cite{apd/sll:igaia05,apd:psr}.

\subsection{Tsallis score}
\label{sec:tsallisscore} The Tsallis score \cite{tsallis:88} is given
by:
\begin{equation}
  \label{eq:tsallisscore}
  S(x,Q) = (\gamma - 1) \int\! d\mu(y)\cdot q(y)^\gamma - \gamma q(x)^{\gamma-1}
  \quad(\gamma>1).
\end{equation}
With minor notational modifications, this is the same as the {\em
  density power score\/} of \textcite{Basu:1998}.

We compute
\begin{equation}
  \label{eq:tsallisent}
  H(P) = -\int\! d\mu(y)\cdot p(y)^\gamma
\end{equation}
and
\begin{equation}
  \label{eq:tsallisdivergence}
  D(P,Q) = \int\! d\mu(y)\cdot p(y)^\gamma + (\gamma-1)\int\! d\mu(y)\cdot q(y)^\gamma - \gamma \int\! d\mu(y)\cdot p(y) q(y)^{\gamma-1}.
\end{equation}
It can be shown that $D(P,Q) > 0$ for $Q\neq P$, demonstrating the
strict propriety of the Tsallis score.

\subsection{Brier score}
\label{sec:brierscore}

Setting $\gamma=2$ in the Tsallis score yields the {\em quadratic\/}
score.  For the special case of a binary sample space ${\cal X} =
\{0,1\}$, an essentially equivalent scoring rule is the {\em Brier\/}
score \cite{brier:50}.  Defining $q := Q(X=1)$ \etc, this has
\begin{eqnarray*}
  S(0, Q) &=& q^2\\
  S(1,Q)  &=& (1-q)^2\\
  H(P) &=& p(1-p)\\
  D(P,Q) &=& (p-q)^2.
\end{eqnarray*}

\subsection{Bregman score}
\label{sec:bregscore} Let $\psi: \reals^+ \rightarrow \reals$ be
convex and differentiable. The associated {\em Bregman score\/} is
given by:
\begin{equation}
  \label{eq:bregscore}
  S(x,Q) = -\psi'\{q(x)\}-\int\!d\mu(y)\cdot\left[\psi\{q(y)\}-q(y)\,\psi'\{q(y)\}\right].
\end{equation}

Then with $p = p(y), q = q(y)$, we get
\begin{eqnarray}
  \label{eq:bregent}
  H(P)& =& -\int \!d\mu(y)\cdot \psi(p),\\
  \label{eq:bregdiv}
  D(P,Q) &=& \int\!d\mu(y)\cdot\left[\psi(p)-\left\{
      \psi(q)+\psi'(q)\,(p-q)\right\} \right].
\end{eqnarray}
By convexity of $\psi$, the integrand of \eqref{bregdiv} is
non-negative, so $S$ is proper (and strictly proper if $\psi$ is
strictly convex).

The log, Tsallis and Brier scores are all special cases of the Bregman
score with, respectively, $\psi(p) = p\ln p$, $\psi(p) = p^\gamma$,
$\psi(p) = (2p^2-1)/4$.

\subsection{Survival score}
\label{sec:survscore}
A variant of the Bregman score, but now applied to the {\em hazard
  function\/} $\lambda_Q(x) := q(x)/\{1-F_Q(x)\}$ (where $F_Q(x) = Q(X
\leq x)$), is useful for scoring a possibly censored survival time
$X$.

Suppose that $X$, non-negative, might be right-censored, at a random
time $C \leq \infty$.  Thus we observe $M = \min\{C,X\}$ and $\Delta =
\bbl(X \leq C)$.  Again let $\psi: \reals^+ \rightarrow \reals$ be
convex and differentiable, and consider the scoring rule:
\begin{displaymath}
  S\{(m,\delta),Q\} = \int_0^m du\cdot\gamma\left\{\lambda_Q(u)\right\}
  - \psi'\left\{\lambda_Q(m)\right\}\,\delta
\end{displaymath}
where \begin{math} \gamma(\lambda) := \lambda\psi'(\lambda) -
  \psi(\lambda).
\end{math}
It can be shown that this is a PSR for the distribution of $X$, even
though observation of $X$ may be subject to an unspecified
non-informative censoring process.

\subsection{Hyv\"arinen score}

Let $\bX$ be a variable taking values in $\cX = \R^k$.  What we term
the {\em Hyv\"arinen score\/} (\textcite{Hyvarinen:2005}; see also
\textcite{Alm-Gid:1993}) is defined by:
\begin{equation}
  \label{eq:hyvarinen}
  S(\bx,Q)=\Delta\ln q(\bx)+\frac{1}{2}\left|\nabla\ln
    q(\bx)\right|^{2}=\frac{\Delta\sqrt{q(\bx)}}{\sqrt{q(\bx)}}
\end{equation}
where $\bnabla$ denotes gradient, and $\Delta$ the Laplacian operator
$\sum_{i=1}^k\partial^2/(\partial x_i)^2$, on $\cX$. With extended
interpretations of these operators, the same expression can be used to
define a proper scoring rule on a general Riemannian space
\cite{apd/sll:igaia05}.

Under conditions that justify ignoring boundary terms when integrating
by parts, we obtain:
\begin{eqnarray*}
  S(P,Q)&=&\frac{1}{2}\int\! d\mu(\by) \cdot \langle\nabla\ln
  q(\by)-2\nabla\ln p(\by), \nabla\ln q(\by)\rangle\\
  H(P)&=&-\frac{1}{2}\int\! d\mu(\by) \cdot\left|\nabla\ln   p(\by)\right|^{2}\\
  D(P,Q)&=&\frac{1}{2}\int\! d\mu(\by) \cdot\left|\nabla\ln p(\by)-\nabla\ln q(\by)\right|^{2}.
\end{eqnarray*}
Since $D(P,Q)>0$ for $Q \neq P$, the Hyv\"arinen score is strictly
proper.  This score also has other important properties that we
highlight in \secref{locality} below.

\subsection{Composite score}
\label{sec:complik}

Consider a model for a multidimensional variable $\bX$.  Let
$\{\bX_k\}$ be a collection of marginal and/or conditional variables,
and let $S_k$ be a PSR for $\bX_k$.  Then we can construct a {\em
  composite score\/} for $\bX$ as
\begin{equation}
  \label{eq:sumscore} S(\bx, Q) = \sum_k S_k(\bx_k, Q_k)
\end{equation}
where $\bX_k \sim Q_k$ when $\bX \sim Q$.  It is easy to see that this
defines a PSR.  It will be strictly proper when every $S_k$ is
strictly proper and the joint distribution for $\bX$ is determined by
the collection of distributions for the $\{\bX_k\}$.

The form \eqref{sumscore} localises the problem to the $\{\bX_k\}$,
which can often simplify computation.  In the special case that each
$S_k$ is the log score, \eqref{sumscore} defines a negative log {\em
  composite likelihood\/} (see \eg\ \textcite{statsinica}).  We can
thus treat composite likelihood in its own right, as supplying a
proper scoring rule, rather than as an approximation (generally poor)
to true likelihood.  Most of the extensive theory and many
applications of composite likelihood apply virtually unchanged to the
more general composite score \eqref{sumscore}.

\subsection{Pseudo score}

A {\em pseudo score\/} is a special case of a composite score.

Consider a spatial process $\bX = (X_v : v \in V)$, where $V$ is a set
of lattice sites.  For a joint distribution $Q$ for $\bX$, let $Q_v$
be the conditional distribution for $X_v$, given the values of
$\bX_{\setminus v}$, the variables at all other sites.  Many
interesting spatial processes are defined locally, by specifying
$\{Q_v, v \in V\}$ (which however can not be done arbitrarily, but is
subject to consistency constraints).  In particular, if $Q$ is Markov,
$Q_v$ only depends on the values of $\bX_{\nei{v}}$, the variables at
the sites neighbouring $v$.

We can construct a proper scoring rule as
\begin{equation}
  \label{eq:pseud}
  S(x, Q) = \sum_v S_0(x_v, Q_v),
\end{equation}
where $S_0$ is a PSR for the state at a single site.  This avoids the
need to evaluate the normalising constant of the full joint
distribution Q.


When $S_0$ is the log score, \eqref{pseud} defines the negative log
{\em pseudo-likelihood\/} of \textcite{besag75}.  Again,
pseudo-likelihood has generally been considered as an approximation to
the full likelihood, but can stand in its own right, as a proper
scoring rule.  For binary $X_v$, taking $S_0$ to be the Brier score
forms the basis of the {\em ratio matching\/} method of
\textcite{Hyvarinen:ext}.  Some comparisons can be found in
\textcite{apd/mm:asta}.

\section{Statistical inference}
\label{sec:est}

\subsection{Estimation}
Let $\{P_{\btheta}: \btheta \in \Theta\}$, where $\Theta$ is an open
subset of $\R^p$, be a parametric family of distributions for
$X\in\cX$.  We suppose given a PSR $S$ on $\cX$, and write
$S(x,\btheta)$ for $S(x,P_\btheta)$, and $s(x,\btheta)$ for its
gradient vector (assumed henceforth to exist) with respect to
$\btheta$:
\begin{eqnarray*}
  s(x,\btheta) &:=& \nabla_\btheta S(x,\btheta)\\
  &=& \left(\frac{\partial S(x, \btheta)}{\partial \theta_j}: j=1,\ldots,p\right).
\end{eqnarray*}

Let $(x_{1},\ldots,x_{n})$ be a random sample from $P_\btheta$,
and denote by $\widehat P$ the empirical distribution of the sample,
which puts mass $1/n$ at each of its (possibly repeated) values.  We
might estimate $\btheta$ by that value minimising $D(\widehat P,
P_\btheta)$, where $D$ is the discrepancy associated with $S$.
Equivalently, since $D(\widehat P, P_\btheta) = S(\widehat P,
P_\btheta) - S(\widehat P, \widehat P)$, we minimise $n S(\widehat P,
P_\btheta)$, which is just the total {\em empirical score\/},
$\sum_{i=1}^n S(x_i, \btheta)$.  That is, our estimate is
$$\widehat{\btheta}_{S}=\arg\min_{\btheta}\sum_{i=1}^{n}S(x_i,\btheta),$$
which (if it exists, which we here assume) will be a root of the {\em
  score equation\/}:
\begin{equation}
  \label{eq:1}
  s(\btheta) :=\sum_{i=1}^{n}s(x_i,\btheta) = \bo.
\end{equation}
We call $\widehat{\btheta}_S$ the {\em minimum score estimator\/} of
$\btheta$.  Note that when $S$ is the log score the score equation is
just the (negative of) the likelihood equation, and the minimum score
estimate is just the maximum likelihood estimate.

Generalising a familiar property of the likelihood equation, the
following theorem \cite{apd/sll:igaia05} shows that, for any proper
scoring rule, and any family of distributions, the score equation
\eqref{1} constructed as above will yield an unbiased estimating
equation:

\begin{theorem}
  \label{teo:1} $$\E_{\btheta}s(X,\btheta)=\bo.$$
\end{theorem}

As a consequence of this theorem we have that equation (\ref{eq:1})
delivers an M-estimator \cite{Hub-Roch:2009}.  We can thus apply
standard results on unbiased estimating equations to describe the
properties of the minimum score estimator $\widehat{\btheta}_{S}$. In
particular, this estimator is consistent in repeated \iid\ sampling.

Define
\begin{eqnarray}
  \label{eq:J}
  J(\btheta) &=& \E_{\btheta}\left\{s(X,\btheta)s(X,\btheta)^T\right\},\\
  \label{eq:H}
  K(\btheta) &=& \E_{\btheta}\left\{ \nabla_\btheta s(X,\btheta)^T\right\},
\end{eqnarray}
with entries
\begin{eqnarray}
  \label{eq:Jab}
  J(\btheta)_{ab} &=& \E_{\btheta}\left\{\frac{\partial S(X,\btheta)}{\partial \theta_a}
    \frac{\partial S(X,\btheta)}{\partial \theta_b}\right\},\\
  \label{eq:Hab}
  K(\btheta)_{ab} &=& \E_{\btheta} \left\{\frac{\partial^2 S(X,\btheta)}{\partial \theta_a\partial \theta_b}\right\},
\end{eqnarray}
and introduce the {\em Godambe information matrix\/}:
$$G(\btheta) :=K(\btheta)J(\btheta)^{-1}K(\btheta).$$
Then under regularity conditions on the model \cite{b-n:1994}, our
estimator is asymptotically normal, with asymptotic covariance matrix
given by the inverse Godambe information matrix:
$$\widehat{\btheta}_{S} \approx \norm(\btheta, \left\{nG(\btheta)\right\}^{-1})$$
when $X_{1},X_{2},\ldots,X_{n}$ are \iid\ as $P_\btheta$.

\subsection{Robust estimation}

The influence function (\IF) of the estimator $\widehat{\btheta}_{S}$,
the solution of the unbiased estimating equation \eqref{1} deriving
from the PSR $S$, measures the effect on the estimator
of adding an infinitesimally small amount of contamination at the
point $x$.  It is given by \cite{Hub-Roch:2009}:
\begin{equation}
  \label{eq:IF}
  \IF_S(x;\btheta) = K(\btheta)^{-1} s(x,\btheta).
\end{equation}

Of particular importance is the supremum of the influence function
over all $x$, a measure of the worst-case influence on
$\widehat{\btheta}_{S}$ of contamination in the data.  For a robust
estimator, this supremum should be finite, \ie, for fixed $\btheta$,
$\IF_S(x;\btheta)$ should be bounded --- this property defines {\em
  B-robustness\/}.  From \eqref{IF} we see that that this will obtain
if and only if the function $s(x;\btheta)$ is bounded in $x$ for each
$\btheta$.

The influence function can also be used to evaluate the asymptotic
variance, $\{nG(\btheta)\}^{-1}$, of $\widehat{\btheta}_{S}$:
$$G(\btheta)^{-1} = \E_\btheta \left\{ \IF_S(X; \btheta) \,
  \IF_S(X;\btheta)^T \right\} \ .$$

\subsubsection{Example: location model}

Suppose $\cX=\Theta=\R$, and the Lebesgue density $p_\theta(\cdot)$ of
$P_\theta$ is given by
\begin{displaymath}
  p_\theta(x) = f(x - \theta),
\end{displaymath}
where the function $f$ is positive and differentiable on $\R$.  We
consider estimation based on the Bregman score \eqref{bregscore} for
given function $\psi$.  We find
\begin{equation}
  \label{eq:rob}
  s(x, \theta)= \psi''\left\{f(u)\right\}f'(u)
\end{equation}
where $u = x-\theta$.  In particular, for the Tsallis score, with
$\psi(t) = t^\gamma$, the necessary and sufficient condition for
B-robustness is that $f(u)^{\gamma-2} f'(u)$ be a bounded function of
$u$ \cite{Basu:1998}.  This condition is satisfied for the normal
location model.

Expression \eqref{rob}, together with the fact that boundedness of
$f'$ implies boundedness of $f$ (see \textcite{DMV:2014}), suggest the
following sufficient conditions for B-robustness:

\begin{cond}
  \label{cond:sc}
  \quad\vspace{-4ex}\\
  \begin{enumerate}
  \item\label{it:f} $f'(u)$ is bounded.
  \item\label{it:psi} $\psi''(t)$ is bounded on $(0,M]$ for any
    $M\in(0,\infty)$.
  \end{enumerate}
\end{cond}
\condref{sc}.\itref{f} holds, for example, for $f$ the normal, the
logistic, the Cauchy or the extreme value distribution.  In typical
cases, \condref{sc}.\itref{psi} will hold so long as $\psi''(0) <
\infty$.

The Brier score, with $\psi(t)= (2t^2-1)/4$, satisfies
\condref{sc}.\itref{psi}: indeed, $\psi''(t) \equiv 1$ is bounded on
the whole of $(0,\infty)$.  For $\gamma> 2$ the Tsallis score
satisfies \condref{sc}.\itref{psi} with $\psi''(0) = 0$.  However for
the log score, having $\psi(t) \equiv t\ln(t)$, $\psi''(t) \equiv 1/t$
is not bounded at $0$, so this particular Bregman scoring rule
violates \condref{sc}.\itref{psi}.  This is reflected in the fact that
the maximum likelihood estimator is typically not B-robust.




\section{Evading the normalising constant}
\label{sec:evading}

When we use the log score, \eqref{1} is just the likelihood equation,
and we obtain the maximum likelihood estimator.

Often we will know the density $p_\theta$ only up to a multiplier:
$$p(x \mid \theta)\propto  f(x \mid \theta)$$
where the omitted normalising constant, $Z(\theta) := \int\!
d\mu(y)\cdot f(y \mid \theta)$, may depend on $\theta$, but not on
$x$.  In this case to solve \eqref{1} we generally need to be able to
compute and differentiate $Z(\theta)$, but often this cannot be done
explicitly.  The identical problem affects estimates based on Bregman
scores and many others.

One solution to this problem proposed in the literature is to use a
composite likelihood approach, which will often avoid the requirement
to evaluate and manipulate $Z(\theta)$.  We will see below that an
alternative escape route is possible by using a suitable {\em local\/}
PSR.

\subsection{Locality}
\label{sec:locality}

To evaluate the log score we only need to know the value of Your
forecast density function, $q(\cdot)$, at the value $x$ of $X$ that
Nature in fact produces.  It is thus termed a {\em strictly local\/}
proper scoring rule.  It can be shown that this property essentially
characterises the log score. However, we can slightly weaken the
locality requirement to admit further PSRs.  For the case of a sample
space that is a real interval, we ask that $S(x, Q)$ should depend on
$q(\cdot)$ only through its value and the value of a finite number of
its derivatives at $x$.  \textcite{mfp/apd/sll:plsr} have
characterised all such local PSRs as a linear combination of the log
score and what they term a {\em key local\/} scoring rule, having the
form
\begin{equation}
  \label{eq:keyc}
  S(x, Q) = \sum_{k=0}^t  (-1)^{k}  \frac{{\rm d}^k }{{\rm d} x^k}
  \phi_{[k]}\left\{x, q(x), q'(x), \ldots, q^{(t)}(x)\right\},
\end{equation}
where $\phi(x, q_0,\ldots,q_t)$ is $1$-homogeneous (\ie, $\phi(x,
\lambda q_0,\ldots,\lambda q_t) \equiv \lambda \phi(x,
q_0,\ldots,q_t)$ for all $\lambda>0$) and concave in
$(q_0,\ldots,q_t)$ for each fixed $x$, and $ \phi_{[k]}$ denotes
$\partial\phi/\partial q_k$.  Some multivariate extensions are
considered by \textcite{mfp:wsc}.

The simplest key local scoring rule is the Hyv\"arinen score, given by
\eqref{hyvarinen} with $k=1$, which arises on taking $\phi = -
q_1^2/q_0$ in \eqref{keyc}.

An important property of every key local scoring rule is {\em
  homogeneity\/}: it is unchanged if $q(\cdot)$ is scaled by a
positive constant.  In particular, $S(x,Q)$ can be computed without
knowledge of the normalising constant of the distribution $Q$.  Thus
if the main computational challenge is to compute this normalizing
constant, it can be tackled by applying a homogeneous scoring rule to
the full joint distribution.

\subsection{Example: Markov process}
Consider the following Gaussian dispersion model for a vector $\bY$
taking values in $\R^N$:
$${\bY} \sim \norm(\bo, \Phi^{-1})$$
with
\begin{equation}
  \label{eq:tridiag}
  \Phi\,\,\,(N \times N)   = \left(
    \begin{array}[c]{cccccc}
      \alpha & \beta & 0 & 0 & \cdots & 0\\
      \beta & \alpha & \beta & 0 & \cdots & 0\\
      0 & \beta & \alpha & \beta& \cdots & 0\\
      \vdots & \vdots & \ddots & \ddots & \ddots & \vdots \\
      0 & 0 & \cdots & 0  & \beta  &\alpha
    \end{array}
  \right)
\end{equation}
where, to ensures that $\Phi$ is positive definite, we take the
parameter space to be
$$\Omega = \{(\alpha,\beta): \alpha > 2|\beta|\}.$$
Note that $\alpha^{-1}$ is the residual variance of each $Y_i$, given
its neighbours.  This model describes a Gaussian time series that is
Markov and approximately stationary.

The determinant of $\Phi$ is
\begin{equation}
  \label{eq:dett}
  \det\left(\Phi\right) = \beta^N \,\frac{\rho^{N+1} - \rho ^{-(N+1)}}{\rho - \rho^{-1}}
\end{equation}
where $\rho$ is determined by
\begin{equation}
  \label{eq:rhoo}
  \rho + \rho^{-1} = \alpha/\beta.
\end{equation}


For an observed data-sequence $\bY=\by$, the likelihood is
proportional to
$$\det\left(\Phi\right)^{\half} \exp \left(-\half \by\transp  \Phi \by\right)$$
with $\det(\Phi)$ given by \eqref{dett} and \eqref{rhoo}.  This will
be hard to maximise directly.

The Hyv\"arinen score \eqref{hyvarinen} eliminates the problematic
normalising constant, and yields a simple quadratic:
\begin{equation}
  \label{eq:quad}
  S(\alpha,\beta) =  - N \alpha + \half \sum_{i=1}^N (\alpha y_{i} +\beta z_i)^2
\end{equation}
where $z_i := y_{i-1} + y_{i+1}$ (taking $y_{-1} = y_{N+1} = 0$).  So
it is easy to minimise directly.  (Note however that the unconstrained
minimum might not belong to $\Omega$, in which case the minimum score
estimate does not exist).

Defining $\lambda = -\beta/\alpha$, \eqref{quad} is
\begin{equation}
  \label{eq:lam}
  -N\alpha + \half\alpha^2 \sum_{i=1}^N(y_i - \lambda z_i)^2.
\end{equation}

The unconstrained minimum is given by
\begin{eqnarray}
  \label{eq:lamhyv}
  \widehat \lambda &=& \frac{c_{yz}}{c_{zz}}\\
  \label{eq:alphhyv}
  \widehat \alpha^{-1} &=& \frac {c_{yy.z}} N
\end{eqnarray}
(and then $\widehat \beta = - \widehat \alpha \widehat \lambda$),
where $c_{yz} := \sum_{i=1}^N y_iz_i$ \etc, and $c_{yy.z} := c_{yy} -
(c_{yz})^2/c_{zz}$.  These will be the minimum score estimates so long
as they lie in $\Omega$, which holds when $c_{yz}^2 < c_{zz}^2/4$.

Alternatively we can apply pseudo-likelihood.  The full conditionals
are given by
$$Y_i | ( \bY_{-i}= \by_{-i}) \sim \norm\left(\lambda z_i, \alpha^{-1}\right)
\quad(i=1, \ldots, N),$$ and the log pseudo-likelihood is thus, up to
a constant:
\begin{equation}
  \label{eq:PL}
  \half N \log \alpha - \half \alpha \sum_{i=1}^N(y_i - \lambda z_i)^2.
\end{equation}
Maximising this gives the same estimates as for the Hyv\"arinen score.

\subsubsection{Multiple observations}

Now suppose we have $\nu$ independent vectors $\bY_1, \ldots,
\bY_\nu$, all distributed as $\norm(\bo,\Phi^{-1})$.  We could form an
estimating equation by summing those derived for the individual
vectors, using either the Hyv\"arinen or the log pseudo-likelihood
score.  This leads again to equations \eqref{lamhyv} and
\eqref{alphhyv}, with $c_{yz}$ redefined as
$\sum_{n=1}^\nu\sum_{i=1}^N y_{ni} z_{ni}$, \etc, and $N$ replaced by
$\nu N$ in \eqref{alphhyv}.

However, we note that a sufficient (albeit not minimal sufficient)
statistic in this problem is the sum-of-squares-and products matrix $S
= \sum_{n=1}^\nu \bY_n \bY_n\transp$, which has a Wishart
distribution: $S\sim W_N(\nu; \Phi^{-1})$.  And the above estimates
are not a function of $S$.  To construct more efficient estimators, we
take $S$ as our basic observable.

It is not clear how pseudo-likelihood could be applied to this
problem.  However, we can still apply a homogeneous scoring rule.
Assume $\nu \geq N$, so that the Wishart density exists, and consider
the multivariate Hyv\"arinen score \eqref{hyvarinen} based on
variables $(t_{ij}: 1 \leq i \leq j \leq N)$, where $t_{ii} = s_{ii}$,
and $t_{ij} = s_{ij}/\sqrt{2}$ for $i<j$.  The associated estimate of
$\Phi$ is obtained by minimising
\begin{equation}
  \label{eq:crit2}
  \sum_{i,j} \left\{(\nu-N-1) s^{ij} - \phi_{ij}\right\}^2
\end{equation}
where $s^{ij}$ denotes the $(i,j)$ entry of $S^{-1}$.  If $\Phi$ is
totally unrestricted, this yields the unbiased estimate
$$\widehat \Phi = (\nu-N-1) S^{-1}.$$
Taking $\Phi$ to have the tridiagonal form \eqref{tridiag}, we get
\begin{eqnarray*}
  \widehat \alpha &=& \frac{\nu-N-1}{N}\sum_{i=1}^N{s^{ii}}\\
  \widehat \beta &=& \frac{\nu-N-1}{N-1}\sum_{i=1}^{N-1}s^{i,i+1}
\end{eqnarray*}
(so long as these estimates satisfy $(\widehat \alpha, \widehat \beta)
\in\Omega$).

\section{Bayesian Model Selection}
\label{sec:modsel}

Suppose that the distribution of an observable $X$ is drawn from one
of a discrete collection ${\cal M}$ of competing parametric models,
where under $M$ the density at $X=x$ is $p_M(x \cd \btheta_M)$,
with unknown parameter $\btheta_M\in\R^{d_M}$.

The Bayesian approach requires us to specify, for each $M\in{\cal M}$,
a prior density function $\pi_M(\btheta_M)$ for its parameter
$\btheta_M$.  Of central importance is the marginal density of $X$
under model $M$, given by:
\begin{equation}
  \label{eq:preddens}
  p_M(x) = \int\! d\btheta_M \cdot p_M(x \cd \btheta_M)\, \pi_M(\btheta_M).
\end{equation}

On observing $X=x_0$, the various models can be compared by means
of the {\em marginal likelihood\/} function, $L(M) \propto
p_M(x_0)$.  In particular, the posterior odds in favour of model $M$
as against model $M'$ are obtained on multiplying the corresponding
prior odds by the {\em Bayes factor\/}, $\BF^M_{M'} = L_{M}/L_{M'}$.

The marginal density \eqref{preddens}, and hence the marginal
likelihood, is sensitive to the choice of the prior distribution
$\pi_M$.  Unfortunately this problem is not solved by using so-called
non informative or objective priors.  These priors are typically
improper and specified in the form $\pi_M(\btheta_M) \propto
h_M(\btheta_M)$.  That is to say, $\pi_M(\btheta_M)=c_M
h_M(\btheta_M)$, where $c_M$ is an unspecified constant.  The same
arbitrary scale factor $c_M$ will then appear in the formal expression
\eqref{preddens} for the marginal density.  The Bayes Factor
$\BF^M_{M'}$ computed using such priors is then not defined, since it
will depends on the ratio $c_M/c_{M'}$ of arbitrary positive
constants.  A variety of {\em ad hoc\/} methods have been suggested to
evade this problem (see, among others, \textcite{OH,Ber-Per}).  Here
we propose a different solution, using proper scoring rules.

\subsection{Use of scoring rules}
\label{sec:use}

The negative log marginal likelihood, $-\log p_M(x_0)$, is just the
log score for the predictive distribution $P_M$ at the observation
$x_0$.  We might now consider replacing the log score by some other
proper scoring rule, $S(x,Q)$ and using that to compare the models
\cite{mm/apd:wsc,apd/mm:ba}.  That is, we replace the (negative log)
marginal likelihood function by the {\em (marginal) score function\/},
$\SF(M) = S(x_0, P_M)$.  Correspondingly, the (negative log) Bayes
factor, $-\log \BF^M_{M'}$, is replaced by the {\em score
  difference\/},
\begin{equation}
  \label{eq:sf}
  \SD^M_{M'} := S(x_0,P_M) - S(x_0,P_M').
\end{equation}
We are thus comparing different hypothesised models for $X$ by means
of their associated scores at the observation $x_0$.




\subsection{Homogeneous score}
\label{sec:hom}

In particular, if $S$ is {\em homogeneous\/}, then $\SF$ and $\SD$
will be insensitive to the arbitrary choice of scale factor in an
improper prior density, and will deliver a well-defined value --- so
long only as $p_M$ given by \eqref{preddens} is finite at $\bx_0$ (but
need not be integrable over $\bx$).  This is just the condition for
having a proper posterior density $\pi(\theta_M \cd \bx_0)$.  There is
then no impediment to adopting improper non-informative priors, so
obtaining an ``objective'' Bayesian model comparison criterion.

For simplicity we here just consider the use of the Hyv\"arinen score
$S_H$ of \eqref{hyvarinen}.  For the general case of multivariate
$\bX$, we find
\begin{displaymath}
  S_H(\bx,P_M) = \E \left\{\left.S_H\left(\bx, P_{\btheta_M}\right)\right| \bX=\bx \right\}
  +\sum_i \var \left\{\left.\frac{\partial\ln p_M(\bx \mid \btheta_M)}{\partial
        x_i}\right| \bX=\bx \right\}
\end{displaymath}
where expectation and variance are taken under the posterior
distribution of $\btheta_M$ given $\bX = \bx$ in model $M$.  This
score is thus well-defined so long as the posterior is proper (even
though the prior may not be), and the required posterior expectation
and variance exist.

 \begin{ex}
   Suppose the statistical model is an exponential family with natural
   statistic $\bT = \bt(\bX)$:
   \begin{equation}
     \label{eq:expgen}
     p(\bx \mid \btheta) = \exp\left\{a(\bx) + b(\btheta) + \btheta\transp \bt(\bx)\right\}.
   \end{equation}
   Define $\bmu \equiv \bmu(\bx)$, $\Sigma \equiv \Sigma(\bx)$ to be
   the posterior mean-vector and dispersion matrix of $\btheta$, given
   $\bX=\bx$.  Then the multivariate Hyv\"arinen score is given by
   \begin{displaymath}
     S_H(\bx,Q) = 2\Delta a(\bx) + 2\bd\transp\bmu
     +  \left\|\bnabla a(\bx) + J\bmu\right\|^2 + 2\tr J\Sigma J \transp
   \end{displaymath}
   with $\bd \equiv \bd(\bx):= (\Delta t_j)$, $J\equiv J(\bx) :=
   (\partial t_j(\bx)/\partial x_i)$.

 \end{ex}

\begin{ex}
  Consider the following normal linear model for a data-vector $\bY =
  (Y_1, \ldots, Y_N)\transp$:
  \begin{equation}
    \label{eq:lm}
    \bY \sim \norm (X\btheta, \sigma^2 I),
  \end{equation}
  where $X$ $(N \times p)$ is a known design matrix of rank $p$, and
  $\btheta\in\R^p$ is an unknown parameter vector.  We take $\sigma^2$
  as known.

  We give $\btheta$ a normal prior distribution:
  \begin{math}
    \btheta \sim \norm(\bm{m},V).
  \end{math}
  The marginal distribution $Q$ of $\bY$ is then
  \begin{math}
    \bY \sim \norm(X\bm{m}, XVX\transp + \sigma^2 I),
  \end{math}
  with precision matrix
  \begin{eqnarray*}
    \Phi &=&  (XVX\transp + \sigma^2 I)^{-1}\\
    &=& \sigma^{-2}\left\{I - X\left(X\transp X + \sigma^2 V^{-1}\right)^{-1}X\transp\right\}
  \end{eqnarray*}
  on applying equation~(10) of \textcite{dvl/afms}.

  An improper prior can be generated by allowing $V^{-1} \rightarrow
  0$, yielding
  \begin{math}
    \Phi = \sigma^{-2} \Pi,
  \end{math}
  where
  \begin{math}
    \Pi := I - X\left(X\transp X \right)^{-1}X\transp
  \end{math}
  is the projection matrix onto the space of residuals.  Although this
  $\Phi$ is singular, and thus can not arise from any genuine
  dispersion matrix, there is no problem in using it to evaluate the
  Hyv\"arinen score.  We obtain
  \begin{equation}
    \label{eq:normhyvimproper}
    S_H(\by, Q) = \frac 1 {\sigma^4}\{\rss - 2\nu\sigma^2\}
  \end{equation}
  where $\rss$ is the usual residual sum-of-squares, on $\nu:= N-p$
  degrees of freedom.  This is well-defined so long as $\nu > 0$.

  When we are comparing normal linear models all with the same known
  variance $\sigma^2$, \eqref{normhyvimproper} is equivalent to
  $(\rss/\sigma^2) + 2p$, Akaike's \aic for this case --- which is
  known not to deliver consistent model selection.

  An alternative to the multivariate Hyv\"arinen score, which avoids
  this problem, is the {\em prequential\/} Hyv\"arinen score.  This is
  a form of composite score, obtained by cumulating the univariate
  Hyv\"arinen scores for the sequence of predictive distributions of
  each $X_n$, given $(X_1,\ldots,X_{n-1})$.  This yields
  \begin{equation}
    \label{eq:normpreq}
    S_H^N = \sum_{n=p}^N \frac 1 {k_n^{2}\sigma^4} (Z_n^2 -2\sigma^2)\quad(N \geq p)
  \end{equation}
  where $Z_n \sim \norm(0,\sigma^2)$ is the difference between $Y_n$
  and its least-squares predictor based on $(Y_1, \ldots, Y_{n-1})$,
  divided by $k_n$.  Without the term $k_n^2$, \eqref{normpreq} would
  reduce to \eqref{normhyvimproper}, and so be inconsistent.  With it
  (even when $k_n \rightarrow 1$, which will typically be the case),
  the difference between the two expressions tends to infinity, and
  use of $S_H^N$ does indeed deliver consistent model selection.
\end{ex}

\section{Conclusion}
\label{sec:conc} Proper scoring rules, of which there is a very
great variety, supply a valuable and versatile extension to
standard statistical theory based on the likelihood function.
Many of the standard results can be applied, with little
modification, in this more general setting. Homogeneous proper
scoring rules, which do not make any use of normalising constant
of a distribution, prove particularly useful in cases where that
constant is computationally intractable, or even non-existent.  We
have illustrated the application of proper scoring rules for
parameter estimation and Bayesian model selection.  We believe
that there will be many other problems for which they will supply
a valuable additional tool in the statistician's kitbag.

\end{document}